\documentclass[12pt]{article}

\usepackage[T1]{fontenc}    % use new font-encoding (to have bold small caps)
\usepackage{a4wide}
\usepackage{amsmath}        % use AMS packages
\usepackage{amssymb}
\usepackage{amsthm}
\usepackage{paralist}
\usepackage{graphicx}
\usepackage{psfrag}         % write latex output into pictures
\usepackage[off]{crop}
\usepackage{bbm}            % fonts for reals, rationals, ...
\usepackage{url}            % typesetting URLs

\usepackage{calc}           % to do calculations
\usepackage[mathscr]{eucal} % euler fonts for calligraphic letters -> \mathscr{}
\usepackage[small,sc]{caption2}
\usepackage[bookmarks,linkcolor=darkgreen,citecolor=red,pagecolor=darkgreen,filecolor=blue,urlcolor=blue]{hyperref}

\usepackage{ae,aecompl}     % good fonts for pdf output

\usepackage{xspace}      % must be loaded last!

\theoremstyle{plain}
\newtheorem{lemma}{Lemma}[section]

\newtheorem*{maintheorem}{Theorem}

\newtheorem{proposition}[lemma]{Proposition}
\theoremstyle{definition}
\newtheorem{definition}[lemma]{Definition}
\newtheorem{example}[lemma]{Example}
\newtheorem{convention}[lemma]{Convention}  
\newtheorem{observation}[lemma]{Observation}  
\theoremstyle{remark}
\newtheorem*{remark}{Remark} 

\newcommand{\R}{\mathbbm{R}}

\DeclareMathOperator{\conv}{conv}
\DeclareMathOperator{\cone}{cone}

\DeclareMathOperator{\interior}{int}

\DeclareMathOperator{\relint}{relint}
\DeclareMathOperator{\sign}{sign}

\DeclareMathOperator{\vertices}{vert}

\newcommand{\Polar}{\Delta}
\newcommand{\A}{{\mathcal A}}

\newcommand{\calO}{{\mathcal O}}

\renewcommand{\c}{{\boldsymbol{c}}}
\newcommand{\bell}{{\boldsymbol{\ell}}}
\newcommand{\bi}{{\boldsymbol{i}}}
\newcommand{\bj}{{\boldsymbol{j}}}
\newcommand{\bk}{{\boldsymbol{k}}}

\newcommand{\g}{{\boldsymbol{g}}}
\newcommand{\h}{{\boldsymbol{h}}}

\newcommand{\p}{{\boldsymbol{p}}}
\newcommand{\q}{{\boldsymbol{q}}}

\newcommand{\bv}{{\boldsymbol{v}}}
\newcommand{\w}{{\boldsymbol{w}}}
\newcommand{\x}{{\boldsymbol{x}}}

\newcommand{\z}{{\boldsymbol{z}}}
\renewcommand{\O}{{\boldsymbol{0}}}

\graphicspath{{pictures/}}

\begin{document}

\title{On the Monotone Upper Bound Problem}% \\On Small Polytopes}
\author{\textsc{Julian Pfeifle}\thanks{Supported by the European
    Graduate Program \emph{Combinatorics, Geometry, and Computation}
    (GRK 588/2) in Berlin and by the GIF project \emph{Combinatorics
      of Polytopes in Euclidean Spaces} (I-624-35.6/1999)}
\quad and\quad\setcounter{footnote}{6}%
\textsc{G\"unter M.~Ziegler}\thanks{Partially supported by Deutsche
  Forschungs-Gemeinschaft (DFG), FZT86,
  ZI 475/3 and ZI 475/4}\\%Partially supported by DFG}\\
\normalsize TU Berlin, MA 6-2\\
\normalsize D-10623 Berlin, Germany\\
\normalsize \texttt{$\{$pfeifle,ziegler$\}$@math.tu-berlin.de}}
\date{August 19, 2003}
 
\maketitle

\vfill
\begin{abstract}\noindent
  The \emph{Monotone Upper Bound Problem} asks for the maximal number
  $M(d,n)$ of vertices on a strictly-increasing edge-path on a simple
  $d$-polytope with $n$ facets.  More specifically, it asks whether
  the upper bound
  \[ M(d,n)\ \le\ M_{\rm ubt}(d,n) \]
  provided by McMullen's (1970) \emph{Upper Bound Theorem} is tight,
  where $M_{\rm ubt}(d,n)$ is the number of vertices of a
  dual-to-cyclic $d$-polytope with $n$ facets.
%  This problem goes back to the early work of Motzkin (1957)
%  and Klee (1965) on the complexity of the simplex algorithm for linear
%  programming.
  
  It was recently shown that the upper bound $M(d,n)\le M_{\rm
    ubt}(d,n)$ holds with equality for small dimensions ($d\le 4$:
  Pfeifle, 2003) and for small corank ($n\le d+2$: G\"artner et al.,
  2001).  Here we prove that it is \emph{not} tight in general: In
  dimension $d=6$ a polytope with $n=9$ facets can have $ M_{\rm
    ubt}(6,9)=30$ vertices, but not more than $26\le M(6,9)\le29$ vertices
  can lie on a strictly-increasing edge-path.
  
  The proof involves classification
  results about neighborly polytopes, % of small corank, 
  Kalai's (1988) concept of abstract objective functions,
  the Holt-Klee conditions (1998),  %\cite{HoltKlee4}
  explicit enumeration,
  Welzl's (2001) extended Gale diagrams,
  randomized generation of instances, 
  as well as non-realizability proofs via 
  a version of the Farkas lemma.   
\end{abstract}
\eject

\section{Introduction}

In an attempt to understand the worst-case behaviour of the
simplex algorithm for linear programming, 
Motzkin \cite{Mot2} in 1957 considered the maximal number $M_{\rm ubt}(d,n)$
of facets that a $d$-polytope with $n$ vertices could have,
and claimed that the maximum is given by the cyclic $d$-polytopes
$C_d(n)$ with $n$ vertices; by polarity, $M_{\rm ubt}(d,n)$ 
is the maximal number of vertices for a simple $d$-polytope with 
$n$~facets.

Motivated by the same problem, Klee \cite{Klee7} in 1965
asked for the maximal number $M(d,n)$ of vertices that
could lie on a \emph{monotone path} (that is, an edge-path that is
strictly monotone with respect to a linear objective function)
on a $d$-polytope with $n$ facets.

Motzkin's claim was substantiated by McMullen \cite{McM1} in 1970. 
It seems that traditionally McMullen's result,
the \emph{Upper Bound Theorem}, was also 
taken as a solution to Klee's question,
the \emph{Monotone Upper Bound Problem}. However, a priori
it is only clear that for all $n>d\ge2$ one has an inequality
\[
 M(d,n)\ \le\ M_{\rm ubt}(d,n),
\]
but it is not at all clear that equality always holds,
that is, that for all $n>d\ge2$ one can construct a simple
dual-to-neighborly $d$-polytope with $n$ facets that
admits a monotone Hamilton path.
Thus in Ziegler \cite[Problems 3.11* and 8.41*]{Z35} %\cite{Ziegler98}
it was explicitly asked:
\medskip

\noindent
\textbf{The Monotone Upper Bound Problem.}
\emph{How large is $M(d,n)$? \ Does it coincide with $M_{\rm ubt}(d,n)$?}
\medskip

\noindent
The quest for ``bad examples'' for the simplex algorithm
equipped with specified pivot rules has led to exponential
lower bounds for $M(d,n)$.
The most prominent one is $M(d,2d)\ge 2^d$, 
as seen from the famous Klee--Minty cubes \cite{KlMi}. %\cite{Klee-Minty72}.
We refer to Amenta \& Ziegler \cite{AmZi2} %{Amenta-Ziegler98}
for a summary of such lower bounds, formulated in the framework of
``deformed products.'' However, these lower bounds are not tight in
general: For example, for $d=4$ and $n=8$ the Klee--Minty cubes yield
$16\le M(4,8)\le M_{\rm ubt}(4,8)=20$, while indeed $M(4,8)=20$.
Similarly for $d=6$ and $n=9$ one may obtain $24\le M(6,9)\le M_{\rm
  ubt}(6,9)=30$.

However, recently the challenge has been taken up,
and it has been proved that the answer to the second question in the 
Monotone Upper Bound Problem is ``YES,'' that is,
$M(d,n)=M_{\rm ubt}(d,n)$ does hold,
\begin{compactitem}[~$\bullet$~]
\item
for small dimensions, $d\le 4$ 
(Pfeifle \cite{Pfeifle-paths03}), \ and
\item
for small corank, $n-d\le2$ 
(G\"artner, Solymosi, Tschirschnitz, Valtr \& 
Welzl~\cite{gaertner01:_one}). %{Gaertner-etal01}).
\end{compactitem}
In the first case, an interesting aspect is that 
the result cannot be achieved on dual-to-cyclic polytopes,
but more general dual-to-neighborly polytopes are needed.
(These had been missed by Motzkin).
The key to the second result is Welzl's concept of ``extended
Gale diagrams'' that will be crucial for our work as well.

In this paper, we give a detailed analysis of some
cases of corank $n-d=3$.
The main result is that the answer to the Monotone Upper Bound 
Problem is ``NO'' in general: %For $d=6$ and $n=9$ one has
\[ 
26\ \le\ M(6,9)\ <\ M_{\rm ubt}(6,9)\ =\ 30.
\]
Our analysis depends on a combination of a number
of different techniques and results:
\begin{compactitem}[~$\bullet$~]
\item Any polytope with $M_{\rm ubt}(d,n)$ vertices
  is necessarily simplicial dual-to-neighborly. 
  If $n=d+3$ and $d$ is even, then it must be dual-to-cyclic.
%  (This is due to Perles \cite[Sect.~6.2]{Gruenbaum67}, using Gale diagrams.)
\item The graphs $G=G(C_d(n)^\Polar)$ of dual-to-cyclic polytopes are given by
  Gale's evenness criterion. For even $d$, 
  $C_d(n)^\Polar$ has a dihedral symmetry group of order $2n$.
\item Any linear objective function in general position 
  induces an acyclic orientation on $G$,
  which is an ``abstract objective function'' (AOF)
  as introduced by Kalai, and satisfies the Holt-Klee (HK)
  conditions. Moreover, in our case it must induce (and be given by) a
  directed Hamilton path in the graph.
\item 
  The symmetry classes of Hamilton paths that induce HK-AOFs are enumerated 
  by computer.
\item 
  In terms of Welzl's  ``extended Gale diagrams'' %\cite{welzl01:_enter}
  the realizability problem for Hamilton~HK~AOFs is reformulated
  as a problem of $3$-dimensional Euclidean geometry.
\item 
  To prove that some of the Hamilton~HK~AOFs of interest are indeed 
  realizable we use randomized generation methods.
\item 
  To prove non-realizability of AOFs we use a combinatorial technique
  that may be seen as an oriented matroid version (looking at signs only)
  of the Farkas lemma; to obtain short proofs, we have implemented
  automatic search techniques.
\end{compactitem}
Our main findings may be summarized as follows.\enlargethispage{5mm}

\begin{maintheorem}
Let $n=d+3$, $d\ge2$. Then a $d$-polytope with $M_{\rm ubt}(d,n)$ vertices
is necessarily dual-to-neighborly;
if $d$ is even, then it is dual-to-cyclic. 
Hamilton HK AOFs on such polytopes can be classified as follows.
\begin{compactdesc}
\item[\boldmath$d=4$, $n=7$:] 
  There are $7$ equivalence classes of Hamilton~HK~AOFs on 
  the graph of $C_4(7)^\Delta$; exactly $4$ of them are realizable. \\
  In particular, $M(4,7)=M_{\rm ubt}(4,7)=14$. Moreover,
  already for $d=4$ and $n=7$ there are non-realizable
  HK AOFs. (These are the smallest possible parameters.)
\item[\boldmath$d=5$, $n=8$:]
  There are two types of simplicial dual-to-neighborly polytopes;
 % Althshuler \& McMullen \cite{altshulermcmullen73}
  for the dual-to-cyclic one realizable types of Hamilton~HK~AOFs exist.\\
  In particular, $M(5,8)=M_{\rm ubt}(5,8)=20$.
\item[\boldmath$d=6$, $n=9$:]
  There are $6$ equivalence classes of Hamilton~HK~AOFs on 
  the graph of $C_6(9)^\Delta$; none of them are realizable. \\
  In particular, $M(6,9)<M_{\rm ubt}(6,9)=30$.
\end{compactdesc}
\end{maintheorem}
%\begin{compactdesc}
%\item[\boldmath$d=4$, $n=7$:
%  There are $7$ equivalence classes of Hamilton~HK~AOFs on 
%  the graph of $C_4(7)^\Polar$; exactly $4$ of them are realizable. %\\
%  In particular, $M(4,7)=M_{\rm ubt}(4,7)=14$. %
%  Moreover, already for $d=4$ and $n=7$ there are non-realizable HK
%  AOFs. (These are the smallest possible parameters.)
%\\[1mm]
%{\upshape\boldmath$d=5$: }
% % There are two types of simplicial dual-to-neighborly polytopes;
% % Althshuler \& McMullen \cite{altshulermcmullen73}
%  There exist realizable Hamilton~HK~AOFs on the graph of $C_5(8)^\Polar$. %\\
%  In particular, $M(5,8)=M_{\rm ubt}(5,8)=20$.
%\\[1mm]
%{\upshape\boldmath$d=6$: }
%  There are $6$ equivalence classes of Hamilton~HK~AOFs on 
%  the graph of $C_6(9)^\Polar$; none of them are realizable. %\\
%  In particular, $M(6,9)<M_{\rm ubt}(6,9)=30$.
%\end{maintheorem}

\section{The combinatorial model}

%Classification of neighborly polytopes via Gale diagrams
%AOF, HK, Hamilton
%Enumeration of Hamilton~HK~AOFs, use templates
%The $d$-polytopes with at most $d+3$ are polar duals
%of the $d$-polytopes with at most $d+3$ vertices.
%These in turn can be classified in terms of Gale diagrams
%\cite[Chap.~6]{Gruenbaum67} \cite[Chap.~6]{Ziegler98}.

%The classification of neighborly polytopes via Gale diagrams
%We start by describing our models for the combinatorial type of a
%polytope and the orientation induced on its graph by a linear
%objective function.

If a $d$-polytope with $d+3$ vertices is supposed to have the
maximal number $M_{\rm ubt}(d,d+3)$ of facets then it must be
simplicial and neighborly. 
Thus, by polarity, we are looking at simple dual-to-neighborly
$d$-polytopes with $d+3$ facets.

The analysis of such polytopes $P$
is a classical application of Gale diagrams by Perles
\cite[Sect.~6.2]{Gr1}. %\cite{Gruenbaum67}. % \cite[Chap.~6]{Ziegler98}.
%We encode the combinatorial type of a $d$-dimensional
%polytope~$P$ with $n$~facets via the Gale diagram of its polar dual.
%By~\cite[Exercise 7.3.7]{Gruenbaum67}, if $P$~has $M_{\rm ubt}(d,n)$
%vertices (i.e, $P$~is polar-to-neighborly), this diagram is a
%\emph{balanced} configuration of vectors in~$\R^{n-d-1}$: for $n-d=3$,
%we get a plane diagram such that exactly $\lfloor d/2 \rfloor$ resp.\
%$\lceil d/2 \rceil$~vectors lie on each side of the linear span of any
%one of them.
It yields that if
$d\ge2$ is even, then the combinatorial type of~$P$ is uniquely that of
the polar~$C_d(d+3)^\Polar$ of the cyclic $d$-polytope with
$d+3$~vertices. For odd $d\ge3$, more combinatorial types of
simple polytopes exist; for $d=3$ as well as for $d=5$
there is exactly one combinatorial type
in addition to the dual-to-cyclic polytope
(see Altshuler \& McMullen~\cite{altshulermcmullen73}). %McMullen74}.

The following yields our combinatorial model for the
orientations of the graph of~$P$ that may be induced by 
linear objective functions (on some realization of~$P$).
\begin{definition} \label{def:aof}
  On the graph of a simple $d$-polytope~$P$ let
  $\calO$ be an acyclic orientation that has a unique source and sink.
\begin{compactenum}[(a)]
  \item\label{aof:a} 
    $\calO$ is an \emph{AOF~orientation} of~$P$ if
    it has a unique sink in each non-empty face of~$P$.
    In this case $\calO$ also has a unique source in each non-empty face
       (Kalai ~\cite{Ka1}; %{Kalai88g}; 
        Joswig, Kaibel \& K\"orner \cite{JKK02}).
    The orientation is then said to satisfy the \emph{AOF condition}.
    Any linear extension of an AOF~orientation
    is called an \emph{abstract objective function (AOF)} on the vertices
    of~$P$.
  \item\label{aof:b} $\calO$ satisfies the \emph{Holt--Klee conditions}
    (or is an \emph{HK orientation}) if in each $k$-dimensional
    face of~$P$ with $3\le k\le d$ it admits $k$~vertex-disjoint
    directed paths between the unique source and sink.
  \item $\calO$ is an \emph{HK AOF orientation} if it satisfies
    (\ref{aof:a}) and (\ref{aof:b}), and a 
      \emph{Hamilton HK AOF orientation} if it additionally admits a directed Hamilton
    path from source to sink.
  \end{compactenum}
\end{definition}

\noindent
Any linear function in general position (that is, such that no two
vertices have the same value) induces an
AOF~orientation on the graph of~$P$; any such orientation is in
fact an HK~orientation (Holt and Klee~\cite{HoltKlee4}). %{Holt-Klee99}).
The negative of the linear function induces the opposite AOF~orientation.
Any Hamilton AOF orientation induces a \emph{unique}
abstract objective function.

If for some linear function on a $d$-polytope with $d+3$ vertices
there is a monotone path through $M_{\rm ubt}(d,d+3)$ vertices, then the
polytope is simple and dual-to-neighborly, and the linear function
induces a Hamilton HK AOF. So for our problem we have to enumerate
Hamilton HK AOFs on the graphs of dual-to-neighborly $d$-polytopes
with $d+3$ facets, which are in fact dual-to-cyclic in the case of
even dimension.
% We have obtained the following results.

\begin{figure}[hbtp]
    \begin{minipage}[c]{.486\linewidth}
      \centering
      \setlength{\parindent}{0pt}
      $N\!R_1^4$  \hspace{-.45cm}\parbox{.95\linewidth}{
        \includegraphics[width=\linewidth]{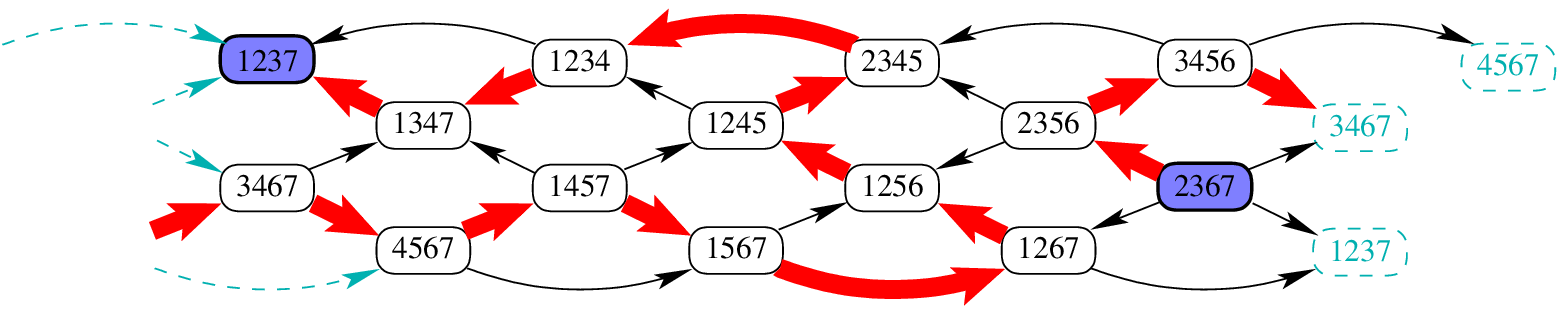}
      }

      \bigskip
      $N\!R_2^4$  \hspace{-.45cm}\parbox{.95\linewidth}{
        \includegraphics[width=\linewidth]{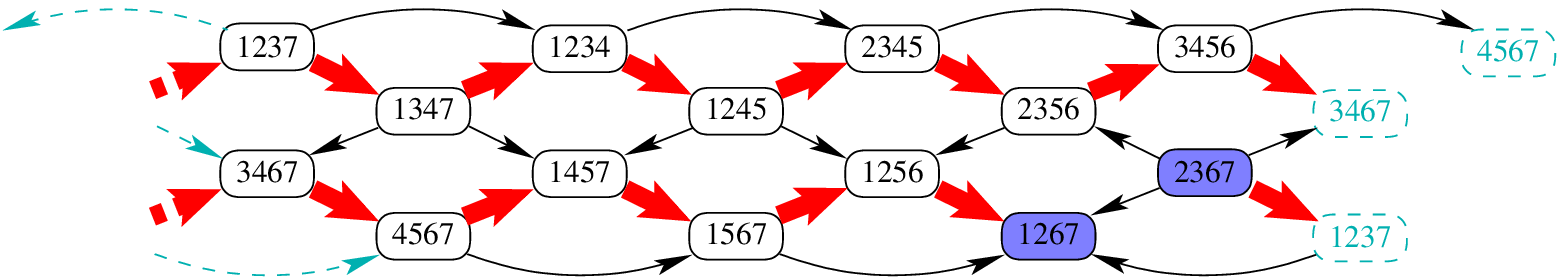}
      }

      \bigskip
      $N\!R_3^4$ \hspace{-.45cm}\parbox{.95\linewidth}{
        \includegraphics[width=\linewidth]{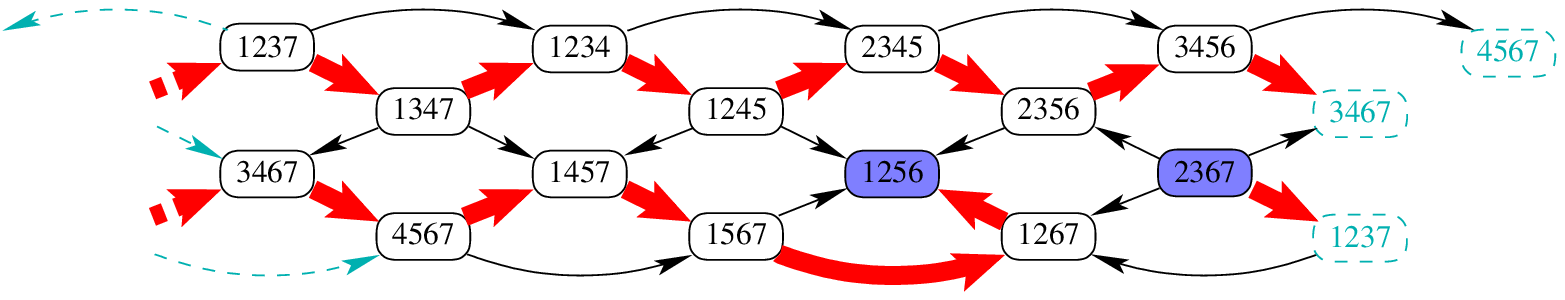}
      }
  \end{minipage}
%  \hfill
  \hspace{-2ex}
  \begin{minipage}[c]{.486\linewidth}
    \centering
      \setlength{\parindent}{0pt}
      \parbox{.95\linewidth}{
        \includegraphics[width=\linewidth]{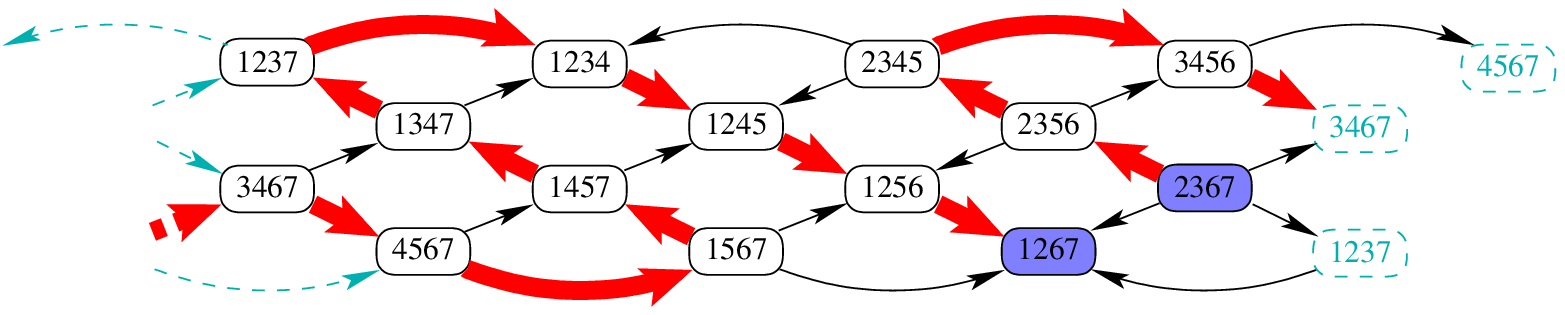}
      }  %\hspace{-.5cm}
      \rlap{$R_1^4$}

%      \smallskip
      \parbox{.95\linewidth}{
        \includegraphics[width=\linewidth]{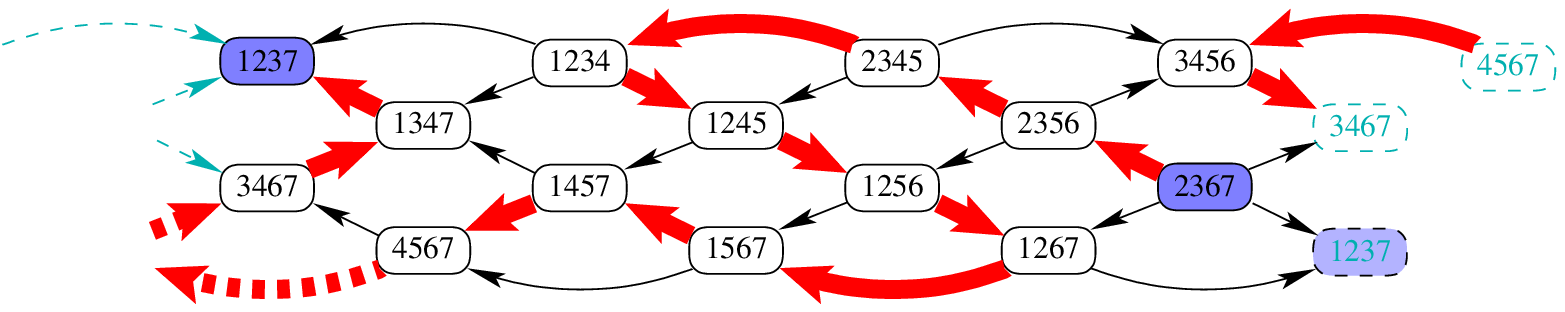}
      }  %\hspace{-.5cm}
      \rlap{$R_2^4$}

%      \smallskip
      \parbox{.95\linewidth}{
        \includegraphics[width=\linewidth]{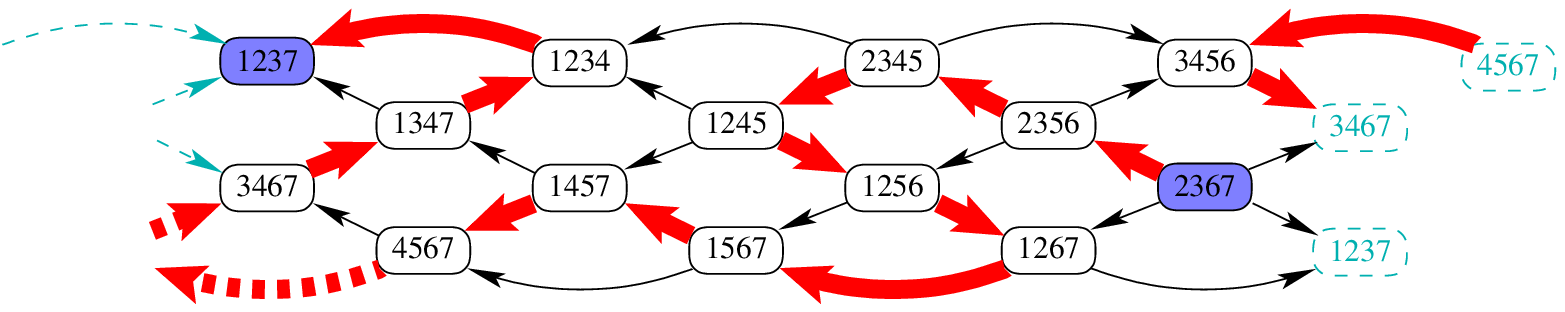}
      }  %\hspace{-.5cm}
      \rlap{$R_3^4$}

%      \smallskip
      \parbox{.95\linewidth}{
        \includegraphics[width=\linewidth]{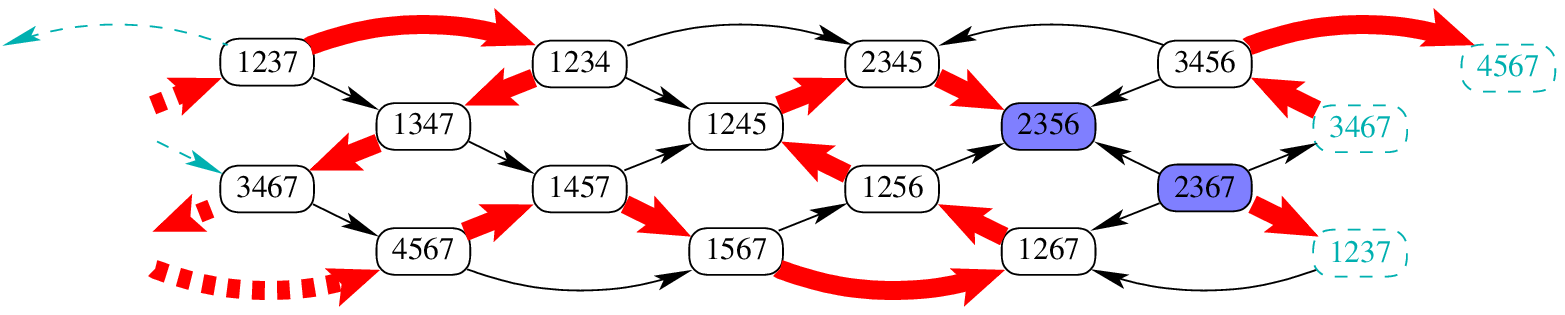}
      }  %\hspace{-.5cm}
      \rlap{$R_4^4$}
  \end{minipage}
  \caption{The Hamilton HK AOFs of
    the graph~$G$ of $C_4(7)^\Polar$. ($G$~embeds into a
    M\"obius strip~\cite{Z76}). %\cite{Haase-Ziegler02}).  
    Each vertex is
    labeled by its set of incident facets, which corresponds to a
    facet of~$C_4(7)$.  The bold arrows yield the monotone Hamilton paths
    from source to sink. An arrow $v\to w$
    means that $w$~is higher than~$v$; 
    so, for example, $N\!R_1^4$ corresponds to\vspace{2mm}\newline
% the following sequence of heights:
    \centerline{$
       2367 < 2356 < 3456 < 3467 < 4567 < 1457 < 1245 < 2345 < 1234 <
       1347 < 1237.
    $}}
  \label{fig.skela1}
\end{figure}

\begin{proposition}\label{prop:enumerate}~
  \begin{compactenum}[\upshape(a)]
  \item The graph of $C_4(7)^\Polar$ admits exactly $7$ equivalence
    classes (with respect to symmetries of $C_4(7)^\Polar$ and global
    orientation reversal) of Hamilton HK~AOFs; 
    they are displayed in Figure~\ref{fig.skela1}.

  \item The polytope $C_5(8)^\Polar$ admits exactly  
%   $2596$ symmetry classes  
    $1298$ equivalence classes of Hamilton  HK~AOFs.
% some of these are realizable.

  \item The polytope $C_6(9)^\Polar$ admits exactly $6$~equivalence
    classes of Hamilton HK AOFs; 
    they are displayed in Figure~\ref{fig:c69}.
% none of these are realizable.
  \end{compactenum}
\end{proposition}

\begin{figure}[htbp]
%  \begin{center}
%\noindent    \begin{minipage}{1.05\linewidth}
%\begin{flushright} 
\upshape\footnotesize%
$N\!R_1^6$: \\ 
\phantom{\quad<} 
458 < 258 < 238 < 278 < 478 < 078 < 058 < 038 < 018 < 014 < 012 < 016 < 
036 < 034 < 345 < \\
\mbox{}\quad < 
234 < 347 < 147 < 127 < 167 < 678 < 367 < 567 < 056 < 456 < 256 < 236 < 
123 < 125 < 145

\medskip
$N\!R_2^6$: \\ 
\phantom{\quad<}
038 < 238 < 123 < 236 < 234 < 034 < 345 < 347 < 478 < 147 < 014 < 018 < 
012 < 016 < 036 < \\
\mbox{}\quad < 
367 < 167 < 678 < 567 < 056 < 256 < 456 < 145 < 458 < 058 < 258 < 125 < 
127 < 278 < 078 

\medskip
$N\!R_3^6$: \\ 
\phantom{\quad<}
038 < 238 < 236 < 036 < 016 < 056 < 256 < 567 < 367 < 167 < 678 < 078 < 
278 < 478 < 147 < \\
\mbox{}\quad < 
127 < 123 < 012 < 125 < 258 < 058 < 458 < 456 < 145 < 345 < 347 < 234 < 
034 < 014 < 018 

\medskip
      $N\!R_4^6$: \\ 
\phantom{\quad<}
038 < 238 < 236 < 036 < 016 < 056 < 256 < 567 < 367 < 167 < 678 < 278 < 
078 < 478 < 147 < \\
\mbox{}\quad < 
127 < 123 < 012 < 125 < 258 < 058 < 458 < 456 < 145 < 345 < 347 < 234 < 
034 < 014 < 018

\medskip
      $N\!R_5^6$: \\ 
\phantom{\quad<}
038 < 058 < 258 < 125 < 256 < 056 < 456 < 458 < 145 < 345 < 034 < 234 < 
347 < 147 < 014 < \\
\mbox{}\quad < 
018 < 012 < 016 < 036 < 236 < 367 < 567 < 167 < 678 < 478 < 078 < 278 < 
127 < 123 < 238 

\medskip
      $N\!R_6^6$: \\ 
\phantom{\quad<}
018 < 058 < 458 < 258 < 125 < 012 < 127 < 278 < 078 < 038 < 238 < 123 < 
234 < 034 < 345 < \\
\mbox{}\quad < 
347 < 478 < 147 < 014 < 145 < 456 < 256 < 236 < 036 < 367 < 678 < 567 < 
167 < 016 < 056 
% \end{flushright}
%\end{minipage}
%\end{center}
\caption{Representatives for the
  six equivalence classes of Hamilton HK AOFs on~$C_6(9)^\Polar$.
  Each vertex~$\p$ is given by the $3$-element set~$N_\p$ of the indices of 
  facets \emph{not} incident to it.}
\label{fig:c69}
\end{figure}

\noindent\emph{Sketch of proof.}
We enumerate the symmetry classes of directed Hamilton paths in the
graph $G$ of one of these polytopes, but prune the search tree
whenever the orientation induced by the partial path fails to satisfy
the AOF or Holt--Klee conditions.

As an additional pruning criterion, we keep a list ${\mathcal L}_F$ of
\emph{all} HK AOF~orientations for each $k$-face~$F$ of~$P$ for
some~$3\le k\le\dim(P)$.  Whenever we try to add a new oriented
edge~$e$ to a partial Hamilton path in~$G$, we check in all lists
$\{{\mathcal L}_F:e\in F\}$ belonging to $k$-faces incident to~$e$
whether there still exists an HK AOF~orientation containing~$e$, and
discard all other orientations of that $k$-face.  
%Without this improvement
%, enumerating all Hamilton AOF
%HK~orientations on~$C_6(9)^\Polar$ is hopeless; even with this idea,
%$C_8(11)^\Polar$ remains far out of reach.

This strategy was implemented in C++ within the \texttt{polymake}
programming environment by Gawrilow \& Joswig 
\cite{GawrilowJoswig,GawrilowJoswig2}; %polymake,Joswig-Gawrilow00});
this produced the results listed above.                   \hfill$\Box$

\section{Extended Gale diagrams}\label{sec:WelzlDiagram}

Welzl's \emph{extended Gale diagram}~\cite{welzl01:_enter,gaertner01:_one}
encodes the values of a linear objective function on a $d$-dimensional
polytope with $n$~facets into an $(n-d)$-dimensional diagram.
%; typically, $n\le d+3$.
For this, we start from a sequence
$(\w_1,\w_2,\dots,\w_n,\g)$ of points in~$\R^d$:
The $\w_i$'s represent the $n$~facet-defining
hyperplanes $\{\x\in\R^d:\w_i^T\!\x=1\}$, $i\in\{1,2,\dots,n\}$, of a
full-dimensional polytope $P\subset\R^d$ with $\O\in\interior P$,
and $\g\in\R^d$ encodes a linear objective function $\g^T\in(\R^d)^*$.

With this interpretation of the input, the extended Gale diagram
produces a sequence $(\w_1^*,\w_2^*,\dots,\w_n^*,\tilde\g^*)$ of
$n+1$~labeled vectors in $\R^{n-d}$ that encodes both the face
lattice of~$P$ and the orientation~$\calO_\g$ of the graph
of~$P$ induced by~$\g^T$. It is calculated as follows:

\begin{compactenum}[(1)]
\item 
  Replace $\g$ by some positive scalar multiple $\tilde\g=c\g$
  such that $\tilde\g^T\!\x<1$ for all $\x\in P$;
  equivalently, $\tilde\g\in\interior P^\Polar$.\\{}
  [This step is optional, and will be modified later. In Welzl's original
  version of extended Gale diagrams it ensures that the ``lifting heights''
  defined below can be made positive.]
%  $P^\Polar=\conv\{\w_1,\w_2,\dots,\w_n\}$ is $P$'s polar dual.
\item Calculate the standard Gale transform $(\w_1^*,\w_2^*,\dots, \w_n^*,
  \tilde\g^*)$ of the point sequence
  $(\w_1,\w_2,\dots,\w_n,\tilde\g)$.
\end{compactenum}

\begin{figure}[ht]
  \centering
  \psfrag{g}{$\tilde\g$} \psfrag{O}{$\O$}
  \includegraphics[height=6cm]{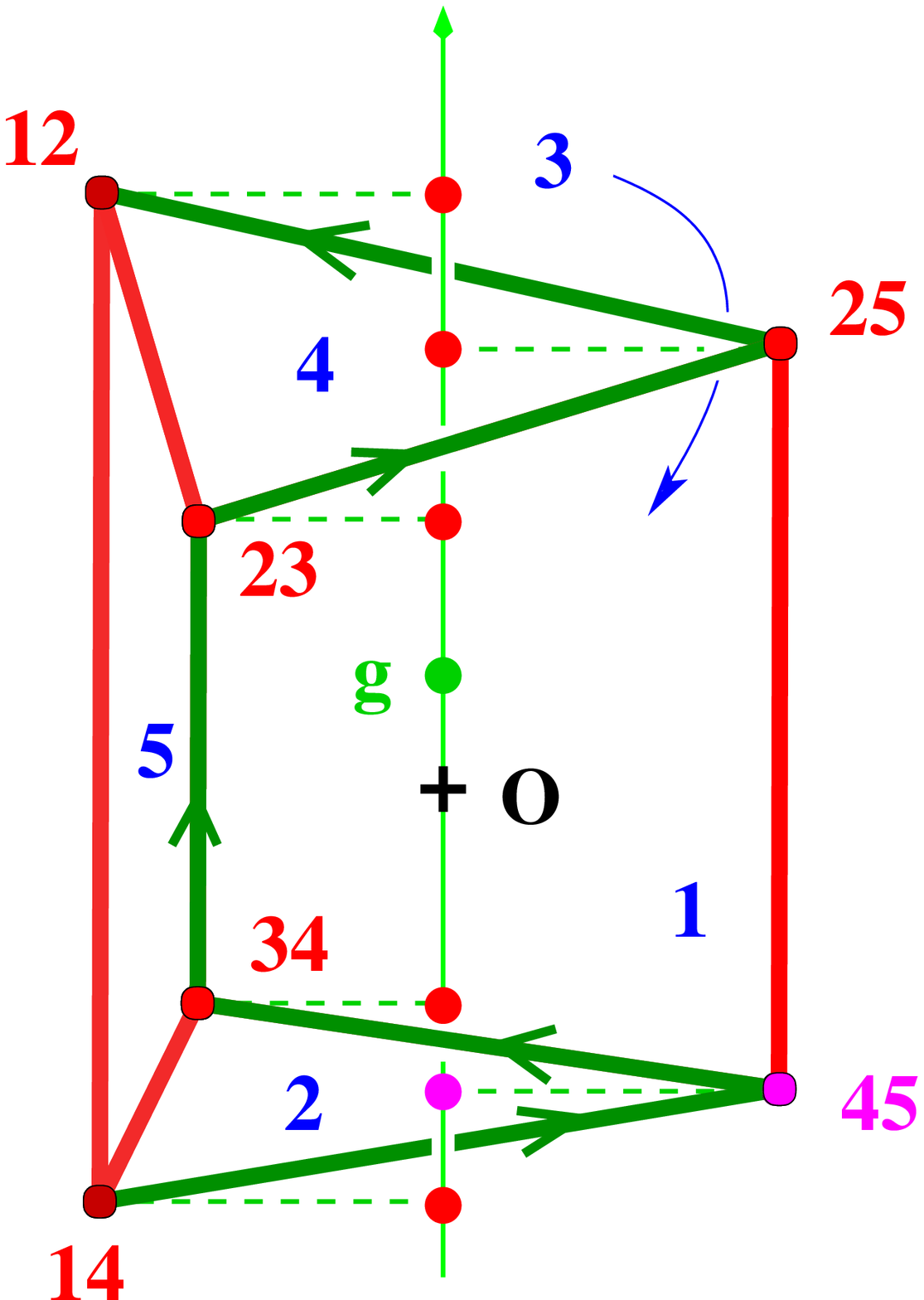}\hfill
  \raisebox{3cm}{$\xrightarrow{\displaystyle\;\Delta\;}$}\hfill
  \raisebox{.5cm}{\includegraphics[height=5cm]{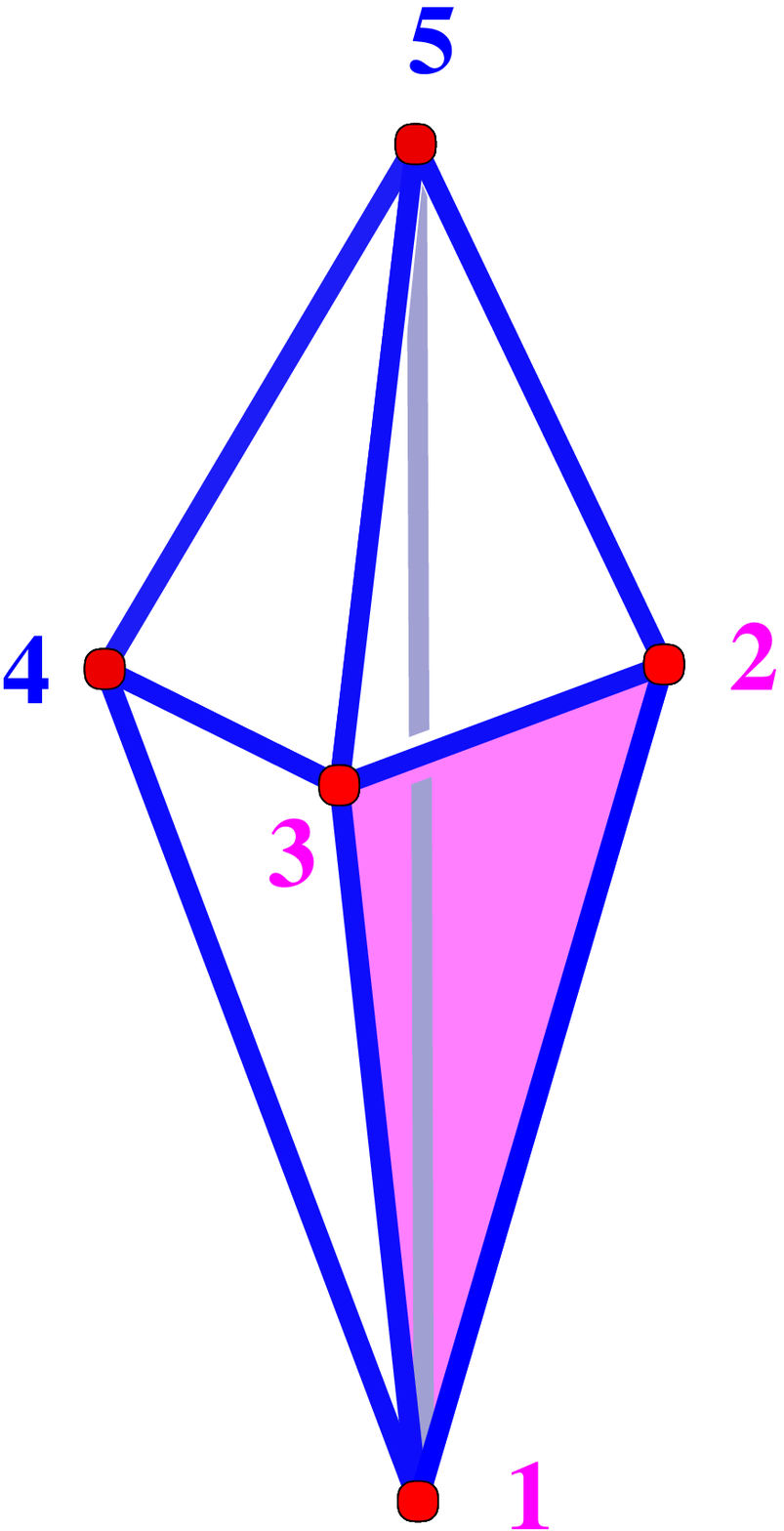}}\hfill
  \raisebox{3cm}{$\xrightarrow{\displaystyle\;\text{Gale}\;}$}\hfill
  \psfrag{g}{$\tilde\g^*$}
  \includegraphics[height=6cm]{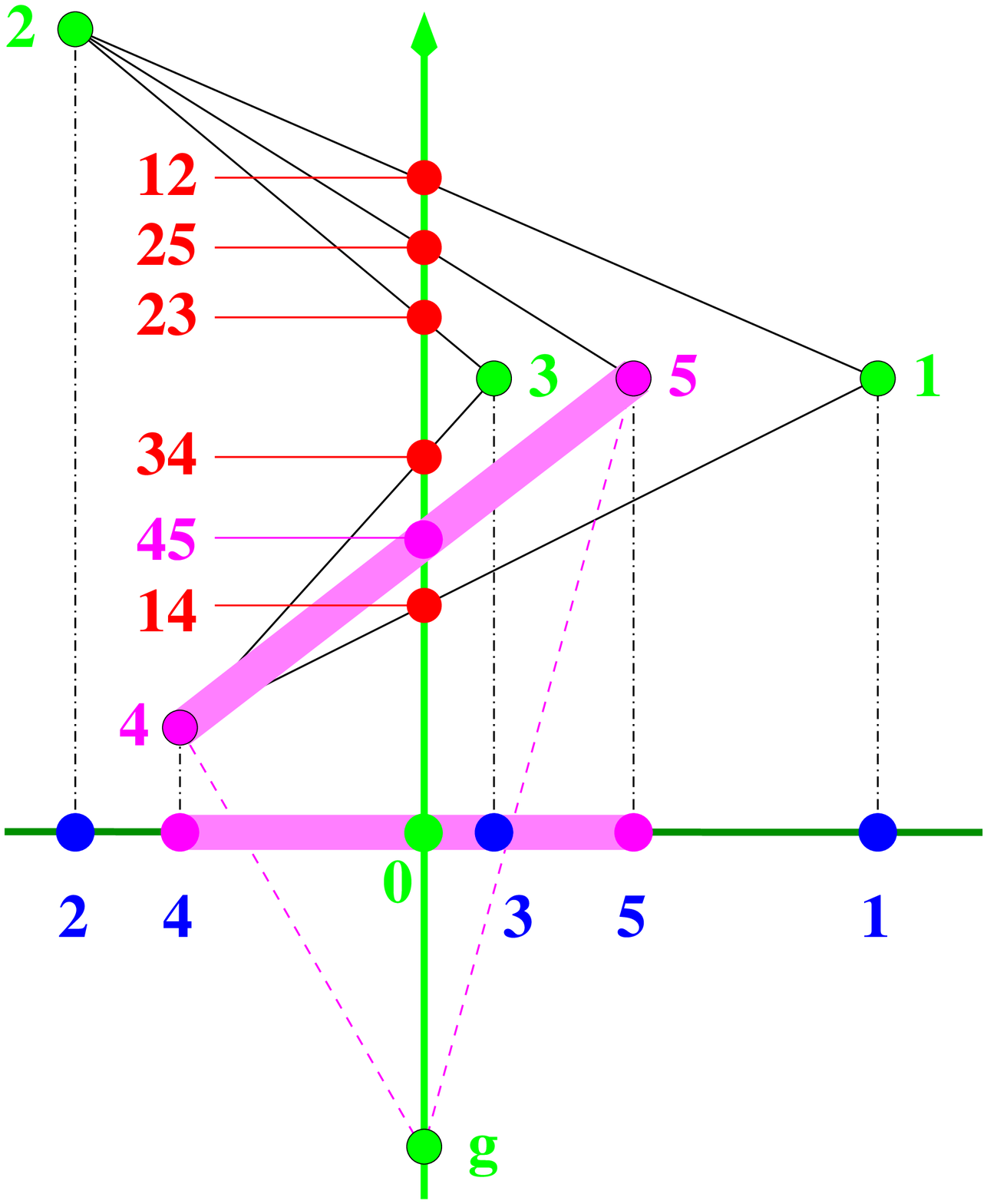}
  \caption{An instance of the extended Gale diagram.\hfill\break
    \emph{Left:} A simple polytope~$P$ whose vertices are labeled with
    the facets they are \emph{not} incident to, and the ordering 
    $14<45<34<23<25<12$ of
    the vertices induced by the linear objective function~$\tilde\g^T$.
    \hfill\break
    \emph{Middle:} The simplicial polar polytope~$P^\Polar$, whose
    vertices are labeled like the corresponding facets of~$P$.\hfill\break
    \emph{Right:} On the base line, a Gale transform of the vertices
    of~$P^\Polar$: Complements of facets of~$P^\Polar$ correspond to
    positive circuits (minimal linear dependencies) of $(\vertices
    P^\Polar)^*$.  Adding~$\tilde\g$ results in a lifting of the Gale
    transform such that the intersection heights for facet complements
    of~$P^\Polar$ encode the ordering of the vertices of~$P$
    by~$\g^T$\!.}
  \label{fig:polar-gale}
\end{figure}

\begin{figure}[ht]
  \centering
  \includegraphics[width=.4\linewidth]{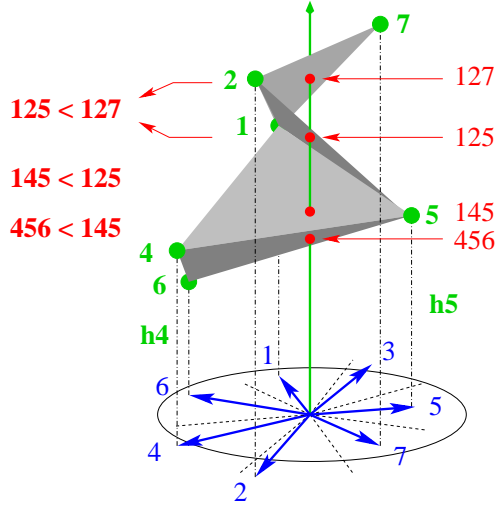}
  \caption{Intersection heights encode values of the objective
    function. Suppose that the objective function $\tilde g$ orders
    four vertices of $C_4(7)$ by $1237 < 2367 < 3467 < 3456$.  Then
    the heights of the intersections between $\R\tilde g^*$~and the
    lifted triangles corresponding to the complements of these labels
    are ordered $z_{456} < z_{145} < z_{125} < z_{127}$. }
  \label{fig:intersection-heights}
\end{figure}

\begin{definition}\label{def:intersection-height}
  Let $P=\{\x\in\R^d:\w_i^T\!\x\le 1,\ 1\le i\le n\}$ be a polytope,
  let $\A_\g=(\w_1,\w_2,\dots,\w_n,\g)\subset\R^d$ be the
  sequence of its facet normal vectors, and let
  $\A^*_\g=(\w_1^*,\w_2^*,\dots,\w_n^*,\tilde\g^*)\subset\R^{n-d}$ 
  be the extended Gale diagram of this sequence,
  whose rows form a basis for the space of affine dependencies
  among the columns of $\A_\g$.

  For every vertex $\p$ of $P$ let $N_\p\subset\{1,2,\dots,n\}$
  index the~$\w_i$ that correspond to the facets of~$P$ that 
  are not incident to~$\p$.  
  The \emph{intersection height}~$z_\p$ of~$\p$ is
  $z_\p=-(\tilde\g^*)^T\!\z_\p$, 
%  the ``negative $\tilde\g^*$-coordinate of $\z_\p$,''
  where $\z_\p=\R\tilde\g^*\cap\conv\{\w_i^*:i\in N_\p\}$ is the
  intersection point of the line~$\R\tilde\g^*$ with the convex hull
  of the~$\w_i^*$'s indexed by~$N_\p$. (See
  Figures~\ref{fig:polar-gale} and~\ref{fig:intersection-heights}.)
%
% The \emph{lifting height} of~$\w_i^*$ is $-(\tilde\g^*)^T\!\w_i^*$.
\end{definition}

\begin{observation}\label{obs:wlog}~
  After a linear transformation we may assume that
  $\tilde\g^*=(0,0,\dots,0,-1)$.  The intersection height $z_\p$ of a
  vertex~$\p$ is then given as the last coordinate of the point where
  the %line $g=\{(0,0,\dots,0,h)^T\!:h\in\R\}$
  $(n-d)$-axis meets the affine plane~$H_\p$ through the points
  $\{\w_i^*:i\in N_\p\}$.
\end{observation}

\begin{proposition}
  \label{prop:intersection-heights} 
  Let $\p,\q$ be vertices of $P$.
  Then $\q$ is higher than $\p$ with respect to the linear objective
  function given by $\g$, 
\[ \g^T\!\p\ <\ \g^T\!\q, \]
if and only the intersection height of $\q$ is larger than that of~$\p$,
\[z_\p\ <\ z_\q.\]
\mbox{~}\vskip-14mm\hfill$\Box$
\end{proposition}

\begin{example}
  Let $P$ be a triangular prism in~$\R^3$ (see
  Figure~\ref{fig:polar-gale}, left) with $n=5$~facets. The
  polar~$P^\Polar$ is the polytope of Figure~\ref{fig:polar-gale}
  (middle) with $5$~vertices, and the Gale transform of~$P^\Polar$
  consists of $5$~points in $\R^{5-3-1}=\R^1$
  (Figure~\ref{fig:polar-gale}, right, base line).  We obtain
  the extended Gale diagram in~$\R^2$ by additionally encoding a
  linear objective function via a level hyperplane that does not
  intersect $P$, which corresponds to a point in the relative interior
  of~$P^\Polar$.  Proposition~\ref{prop:intersection-heights} says
  that in the extended Gale diagram, the value $\g^T\!\p$ of the objective
  function is encoded by the height of the intersection of
  $\R\tilde\g^*$ with the triangle spanned by the points~$\w_i^*$ that
  correspond to facets of~$P$ that do not contain~$\p$.
\end{example}

%\newpage
\section{Finding realizations}

\begin{proposition}\label{prop:realize}~
  \begin{compactenum}[\upshape(a)]
  \item The equivalence classes $R_1^4$--$R_4^4$
    of Hamilton HK AOF~orientations
    of the graph of $C_4(7)^\Polar$ (as given by Figure~\ref{fig.skela1}) 
    are realizable. %\\ 
    In particular,
    $M(4,7)=M_{\rm ubt}(4,7)=14$.
  \item There exist realizable Hamilton HK AOF~orientations
    of the graph of $C_5(8)^\Polar$. \\ In particular,
    $M(5,8)=M_{\rm ubt}(5,8)=20$.
  \item There exist realizations of $C_6(9)^\Polar$ with
    $26$~vertices on a monotone path.
  \end{compactenum}
\end{proposition}

\noindent\emph{Sketch of proof.} 
The realizations were found by the following procedure. For each
polytope $P=C_d(d+3)^\Polar$, randomly generate a Gale
transform~${\mathcal G}(P^\Polar)=(\bv_1^*,\bv_2^*,\dots,
\bv_{d+3}^*)$ of~$P^\Polar$, and for each vertex~$\p$ of~$P$ express
the intersection height~$z_\p$ as a linear function of the
\emph{lifting heights} $h_i$, where $\big(\w_i^*=(\bv_i^*,h_i):1\le
i\le d+3\big)$ is an extended Gale transform of~$P$. Now check whether
the linear program
\[
   z_\p-z_\q\ \le \ -1
   \qquad\text{for all oriented edges $e=(\p,\q)$ in~$\calO$}
\]
in the variables $h_1,h_2,\dots, h_{d+3}$ is feasible, for $\calO$~one
of the Hamilton HK AOF~orientations of~$P$. If so, the polar dual of
the Gale transform of~${\mathcal G}(P^\Polar)$ yields a realization of
the combinatorial type of~$P$, and the lifting heights solving the
linear program yield a linear objective function that induces the
orientation~$\calO$ on this realization.  If not, repeat.  \hfill
$\Box$

%\begin{remark}~
%  \begin{compactitem}[~$\bullet$~]
%  \item There does \emph{not} exist any Hamilton HK~AOF of the graph
%    of~$C_4(8)^\Polar$, so that the non-existence of extremal
%    realizations for this polytope is already proved at the
%    combinatorial level.  However, the other two combinatorial types
%    of polar-to-neighborly $4$-polytopes with $8$~facets~\cite[Section
%    7.2.4]{Gruenbaum67} \emph{do} have such orientations, some of
%    which can even be shown to be realizable by the method of
%    Proposition~\ref{prop:realize}. %Section~\ref{sec:realizable-paths}.
%    In particular, this yields that $M(4,8)=20$.
%  \item The polar-to-cyclic polytopes $C_4(9)^\Polar$ and
%    $C_4(10)^\Polar$ also do not have any Hamilton  HK AOFs
%    (although they \emph{do} admit 165~096 resp.\ 613~040
%    orientations that only fail to satisfy the Holt--Klee condition);
%    we believe this to be true for all polytopes $C_4(n)^\Polar$ for
%    all $n\ge8$, but as yet have no proof.
%%  \item As noted above, there is a second combinatorial type
%%    $N_5(8)$ of a simple polar-to-neighborly polytope with 
%%    $M_{\rm ubt}(5,8)=20$ vertices. 
%%    existence of realizable Hamilton HK~AOFs of the graph
%%    of~$C_5(8)^\Polar$ already shows that $M(5,8)=20$, so that we need
%%    not consider this polytope.
%  \end{compactitem}
%\end{remark}

\section{Proving non-realizability}

%Farkas lemma
Our strategy for proving non-realizability of orientations may be
summarized as follows. For each candidate orientation~$\calO$ of the
graph of a polytope~$P$ (of even dimension $d$, with $d+3$ facets), 
we assume that there is a realization of~$P$ and a linear
objective function~$\g^T$ that induces~$\calO$ on $P$'s graph. Each
oriented edge of~$\calO$ then imposes a linear inequality on the
lifting heights of the extended Gale diagram of~$(P,\g^T)$. For some
orientations~$\calO$, a combinatorial version of the Farkas Lemma
implies that these inequalities are inconsistent, thereby proving the
non-realizability of~$\calO$.

\subsection{Inequalities induced by edges}

We start with some notation for vector configurations in~$\R^2$
and~$\R^3$. The shorthand $[d+3]$ will denote $\{1,2,\dots,d+3\}$.

\begin{convention}\label{conv:vector}
  For $i\in[d+3]$, we write $\bi$ for a vector
  $(x_i,y_i)^T\in\R^2$, and $\bi^\perp$~for the vector $(y_i,-x_i)^T$
  orthogonal to~$\bi$ that is obtained by rotating~$\bi$ in the
  clockwise direction.  With this convention, the following relations
  hold for scalar products:
  \begin{equation*}%\label{eq:ijperp}
    \bi\bj^\perp \;=\; x_iy_j - x_jy_i \;=\; \det(\bi,\bj)
    \;=\; -\det(\bj,\bi) \;=\; -\bj\bi^\perp \;=\; -\bi^\perp\!\bj.
  \end{equation*}
  We further abbreviate
  \begin{equation*}
    ij^\perp \;:=\; \sign(\bi\bj^\perp), \qquad
    \boldsymbol{[ijk]} \;:=\;
    \det\begin{pmatrix} \bi & \bj & \bk \\ 1 & 1 & 1
    \end{pmatrix}, \qquad %\label{eq:ijk} \\{}
    [ijk] \;:=\; \sign(\boldsymbol{[ijk]}). \label{eq:ijk}
  \end{equation*}
%  so that $[ijk]=+$ if and only if $\bi,\bj,\bk$ come in
%  anti-clockwise order around~$0$.
\end{convention}

 \begin{lemma}\label{lem:sign}~
  \begin{compactenum}[\upshape(a)]
  \item If $\bi,\bi+\bj,\bj\in\R^2$ come in anti-clockwise order
    around~$0$, then $ij^\perp=+$.
  \item\label{lem:signs} If in a configuration of four vectors $\bi,
    \bj, \bk, \bell\in\R^2{\setminus}\{0\}$ the vectors $\bi, \bj, \bk$
    are ordered clockwise around~$0$, $\bj\in\relint\cone(\bi,\bk)$,
    $[ijk]=+$, and $\bell\in\relint\cone(-\bi,-\bk)$,
    then $[i\ell j]=[j\ell k]=+$. \hfill $\Box$
  \end{compactenum}
\end{lemma}

\begin{figure}[htb]
  \begin{center}
     \psfrag{i}{$\boldsymbol i$}      \psfrag{j}{$\boldsymbol j$}
     \psfrag{k}{$\boldsymbol k$}      \psfrag{l}{$\boldsymbol\ell$}
     \psfrag{iperp}{$\boldsymbol i^\perp$}
     \psfrag{jperp}{$\boldsymbol j^\perp$}
     \psfrag{hj}{$H_j$}
     \psfrag{h1j}{$H'_j$}
   \includegraphics[height=3.6cm]{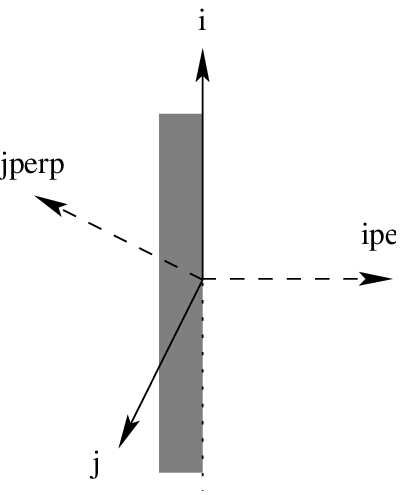} \qquad\qquad\qquad
   \includegraphics[height=3.6cm]{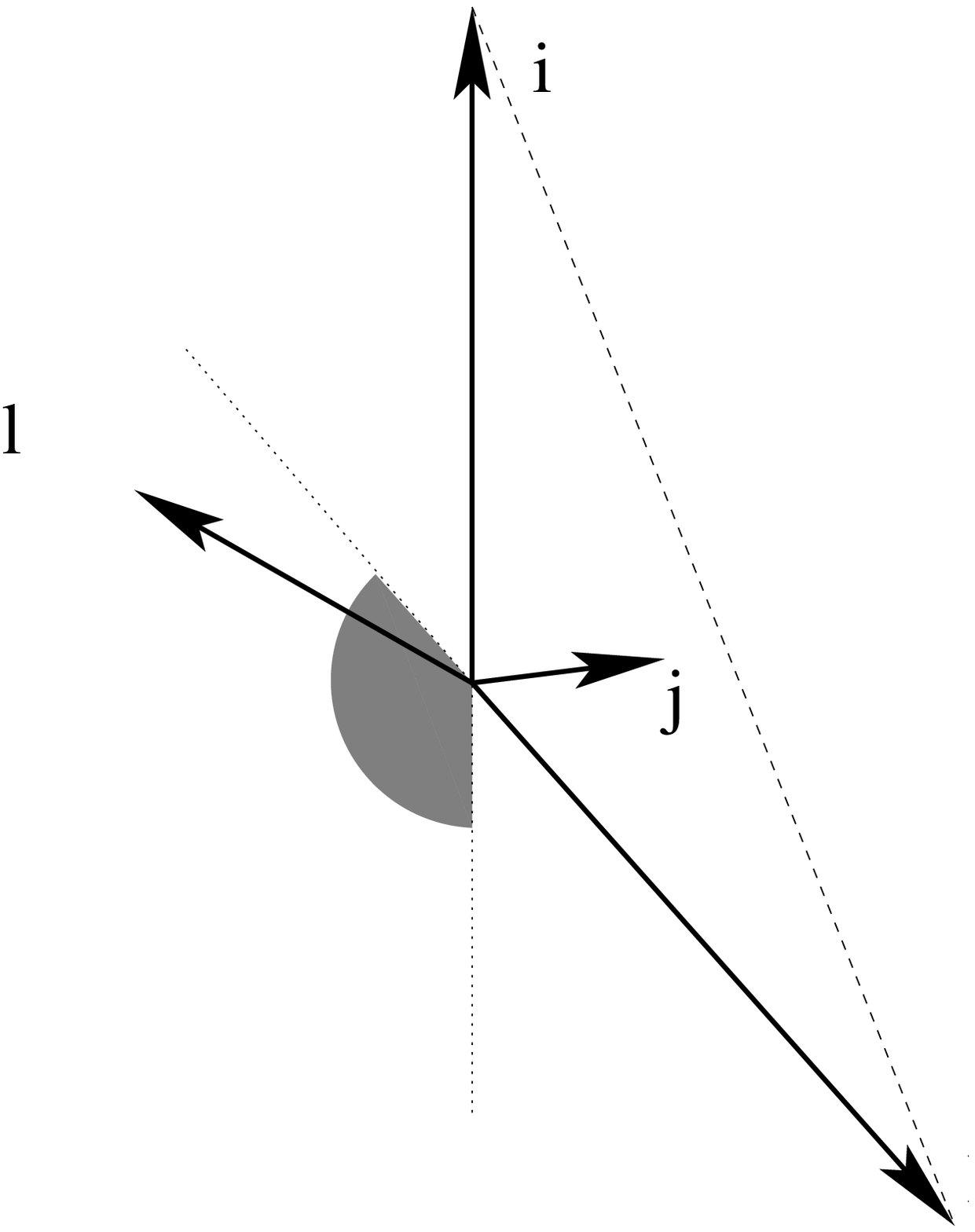}
    \caption{Deducing sign patterns.\hfill\break
    \emph{Left:} If $\bi,\bi+\bj,\bj\in\R^2$ come in
      clockwise order around $0$, then $ij^\perp=+$.\hfill\break 
    \emph{Right:} In this situation, if $[ijk]=+$, 
      then $[i\ell j]=[j\ell k]=+$.}
%      Also, under the conditions of Lemma~\ref{lem:sign}(b), $[ijk]$
%      is the sign of the permutation~$(ijk)$. }
    \label{fig:lemma2}
  \end{center}
\end{figure}

\begin{figure}[ht]
  \centering
  \includegraphics[width=4.4cm]{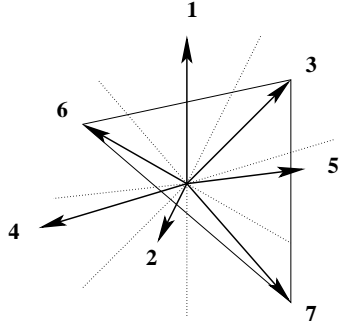}
  \caption{%
    A Gale transform of~$C_4(7)$. The set $N_\p=\{3,6,7\}$ corresponds to
    the vertex~$\p$ of $C_4(7)^\Polar$ not on those facets,
    and $3<6<7$ is an anti-clockwise orientation of the triangle~$367$.}
  \label{fig:convention-gale}
\end{figure}

\begin{convention}\label{conv:gale}
  The vertices of $C_d(d+3)$ are labelled by $[d+3]$ in the natural
  order, so that the facets are given by certain $d$-subsets of
  $[d+3]$ according to Gale's evenness criterion.  The vectors in any
  Gale transform are then labeled so that $\boldsymbol 1, \boldsymbol
  3, \boldsymbol 5, \boldsymbol 7, \dots, \boldsymbol 2, \boldsymbol
  4, \boldsymbol 6, \dots$ come in clockwise order around the origin.
  We identify each facet of $C_d(d+3)$ with the indices of the three
  vertices it misses, so that ordering this index set yields a
  triangle with anti-clockwise orientation that encloses the origin
  (cf.~Figure~\ref{fig:convention-gale}).
  
  Now we polarize.  Correspondingly, we label each vertex $\p$ of
  $C_d(d+3)^\Polar$ by the $3$-element set~$N_\p$ of (indices of)
  the facets it does \emph{not} lie on.
%  In the extended Gale diagram, the vertices of $C_d(d+3)^\Polar$
%  are thus given by oriented affine triangles
%  $\pi_N=\aff\{\w^*_j:j\in N\}$ that intersect the $z$-axis.
%%Thus we also give three element labels also the affine
%%  (hyper)planes $\pi_N=\aff\{\w^*_j:j\in N\}$ in
%%  $\R^{d+3-d}=\R^3$ spanned by points corresponding to complements of
%%  facets of~$C_d(d+3)$ receive the 3-element labels ${\bar
%%    I}=\{1,2,\dots,d+3\}{\setminus} I$. 
\end{convention}

%\item
%  Moreover, for one vertex~$\p$ of~$C_d(d+3)$, with $N_\p=\{i,j,k\}$, we can
%  achieve that the corresponding plane in the 
%  extended Gale diagram is horizontal, that is,
%  that $(\w^*_i)_3=(\w^*_j)_3=(\w^*_k)_3= z_\p$.
%%  $I_\p$~indexes~$F_\p$ (cf.~Definition~\ref{def:intersection-height}).

%\item \textbf{[revise! later??]}
%  It follows that we may assume in any Gale transform of
%  $(C_d(d+3),\g)$ that $h_i=h_j=h_k=0$ for some facet complement
%  $\{i,j,k\}$. The resulting configuration is not a
%  extended Gale diagram of~$C_d(d+3)^\Polar$ and an objective
%  function, but projecting along $\R\tilde\g^*$ does yield a
%  $\Gale$-transform of~$C_d(d+3)$.
%\end{compactenum}
%\end{observation}

\begin{lemma} \label{lem:zijk}
  Let ${\mathcal G}(P)$ be an extended Gale diagram
  of~$P=C_d(d+3)^\Polar$, and let $N_\p=\{i,j,k\}$
  index a vertex~$\p$ of~$P$.
  With the assumptions of Observation~\ref{obs:wlog}
  and Convention~\ref{conv:vector},
  the intersection height $z_\p=z_{\{i,j,k\}}$ %where the line
%  $g=\{(0,0,h)^T\!:h\in\R\}$ meets the affine plane $\pi_{\{i,j,k\}}$
%  through the points $\w_i^*,\,\w_j^*,\,\w_k^*\in\R^3$ (corresponding to
%  the vertex $\{1,2,\dots,d+3\}{\setminus}\{i,j,k\}$ of~$P$)
   is given by
    \begin{equation}\label{eq:z}
      z_{\{i,j,k\}} \;=\; \dfrac{\bi \bj^\perp h_k + \bk \bi^\perp h_j +
        \bj \bk^\perp h_i}{\boldsymbol{[ijk]}}.
    \end{equation}
\end{lemma}

\noindent\emph{Proof.}
Expand the third row of the determinant in the equation
\begin{equation*}
     \left|
       \begin{matrix}
         0 & x_i & x_j & x_k\\
         0 & y_i & y_j & y_k\\
         z_{\{i,j,k\}} & h_i & h_j & h_k\\
         1 & 1 & 1 & 1
       \end{matrix}
       \right| \;=\; 0.
\end{equation*}
\noindent  As a consistency check, note that \eqref{eq:z} is symmetric under any
permutation of the indices.\hfill$\Box$

\medskip

By Proposition~\ref{prop:intersection-heights}, the total ordering of
the vertices $\p$ of $C_d(d+3)^\Polar$ induced by the linear
objective function~$\g^T$ induces a total ordering of the intersection
heights~$z_\p$, that is, of the heights of the intersections of the
affine hyperplanes~$H_\p$ in~$\R^3$ with the $z$-axis. If two vertices
of~$C_d(d+3)^\Polar$ span an edge, then the corresponding facets of
$C_d(d+3)$ share a ridge, which in turn means that the corresponding
triangles have two points~$\w_i^*$, $\w_j^*$ in common.  This permits
us to relate the intersection heights of two adjacent vertices in the
graph of~$C_d(d+3)^\Polar$ in the following way.

\begin{lemma}\label{lem:diff} Suppose that the vertices $\{i,j,k\}$ 
  and $\{i,j,\ell\}$ span an edge of~$C_d(d+3)^\Polar$. Then the
  following relation holds between the corresponding intersection
  heights:
  \begin{equation*}%\label{eq:diff}
    z_{\{i,j,k\}}-z_{\{i,j,\ell\}} \;=\;
\dfrac{(\bi\bj^\perp)\boldsymbol{[jk\ell]}}{\boldsymbol{[ijk][ij\ell]}}h_i+
\dfrac{(\bi\bj^\perp)\boldsymbol{[ki\ell]}}{\boldsymbol{[ijk][ij\ell]}}h_j+
      \dfrac{\bi\bj^\perp}{\boldsymbol{[ijk]}}h_k +
      \dfrac{-\bi\bj^\perp}{\boldsymbol{[ij\ell]}}h_\ell.
  \end{equation*}
  If $[ijk]=[ij\ell]$, then the signs of the coefficients of the $h$'s
  are, in this order,
  \begin{equation*}%    \label{eq:signs}
    (ij^\perp)[jk\ell],\quad (ij^\perp)[ki\ell],\quad
    +,\quad -.
  \end{equation*}
\end{lemma}

\noindent\emph{Proof.} The first statement follows via direct calculation
from~\eqref{eq:z}, using the straightforward identity
\begin{equation*}
  (\bi\bj^\perp)(\bk\bell^\perp) \;=\;
  (\bell\bi^\perp)(\bj\bk^\perp)+ (\bj\bell^\perp)(\bi\bk^\perp).
%\tag*{$\Box$}
\end{equation*}
The second statement is a consequence of Lemma~\ref{lem:sign}
and Convention~\ref{conv:gale}. \hfill $\Box$

\subsection{Contradictions via a combinatorial Farkas Lemma}

We will use a combinatorial version of the following Farkas
Lemma~\cite[Sect.~7.8]{Schr}: %{Schrijver86}:

\begin{lemma}\label{lem:farkas} For any matrix
  $A\in\R^{m\times d}$, exactly one of the following is true:
  \begin{compactitem}[~$\bullet$~]
  \item There exists an $\h\in\R^d$ such that $A\h < \O$.
  \item There exists a $\c\in\R^m$ such that $\c\ge\O$,
    $\c^T\! A=\O$, and $\c\ne\O$.\hfill $\Box$
  \end{compactitem}
\end{lemma}

%\newpage %
Given a $d$-dimensional polytope~$P$ with $d+3$ facets and an
orientation~$\calO$ on $P$'s~graph~$G$, we assume that we have a
realization of~$P$ and a linear objective function~$\g^T$ that
induces~$\calO$ on~$G$. We would like to apply Lemma~\ref{lem:farkas}
to prove the infeasibility of the system $A\h<\O$ of 
$m=\#\textrm{edges of }C_d(d+3)^\Polar=\frac14\binom{d+4}3$ linear
inequalities on the lifting heights $h_1,h_2,\dots,h_{d+3}$ given by
\begin{equation}\label{eq:ineq}
   z_{\{i,j,k\}}-z_{\{i,j,\ell\}} \ <\ 0 \qquad\text{for all oriented
edges }
   (\{i,j,k\},\{i,j,\ell\})\text{ of }\calO.
\end{equation}
However, the only information we have available about~$A$ are sign
patterns of determinants as given by Lemma~\ref{lem:diff}.  Therefore,
to show infeasibility of~\eqref{eq:ineq} we must produce a Farkas
certificate~$\c$ that shows already at the level of signs (``using
only oriented matroid information'') that some positive combination of
the rows of $A$ sums to zero.

\begin{proposition} \label{prop:nr14}
  The orientation
  \begin{eqnarray*}
     N\!R_1^4: \qquad  z_{145} & \!<\!& z_{147} \;<\; z_{127} \;<\;
z_{125} \;<\; z_{123}\;<
    \; z_{236} \;<\; z_{234} \;< \\
     &\!<\!& z_{345} \;<\; z_{347} \;<\; z_{367} \;<\; z_{167} \;<\;
z_{567} \;<\;
       z_{256} \;<\; z_{456}
     \end{eqnarray*}
     of the graph of $C_4(7)^\Polar$ is not realizable.
\end{proposition}

\noindent\emph{Proof.} We abbreviate
`$z_{\{i,j,k\}} < z_{\{i,j,\ell\}}$' by `$ijk < ij\ell$'.

To any extended Gale diagram corresponding to a realization
of~$N\!R_1^4$ we may apply an affine transformation that fixes the
$z$-axis and moves the plane spanned by $\w^*_3,\w^*_4$ and
$\w^*_5$ to the $\R^2$-plane given by $z=0$; that is, we may assume
that $h_3=h_4=h_5=0$. 
This affine transformation does not change the projection
along the $z$-axis, which still yields
the same Gale transform of~$C_4(7)$.
The resulting configuration is 
the extended Gale diagram for $C_4(7)^\Polar$ with 
the objective function $\tilde\g=c\g$ scaled such that the level hyperplane
$\tilde\g^T\!\x=1$ contains $\p=\{3,4,5\}$.
Thus at this point we have modified Step (1) in the construction of
Section~\ref{sec:WelzlDiagram}.

We proceed to write down the sign patterns of
the inequalities $A\h<\O$ for $\h=(h_1,h_2,h_6,h_7)$ implied by
Lemmas~\ref{lem:zijk}~and~\ref{lem:diff}:

\medskip\qquad
\begin{tabular}[b]{c|cccc|*{4}{c@{\;}}}
             & $h_1$ & $h_2$ & $h_6$ & $h_7$
             & $i$ & $j$ & $k$ & $\ell$ \\\hline
  $567<256$: &   0   &  $-$  & $-[257]$ &  $+$  & 5&6&7&2 \\
  $234<345$: &   0   &  $+$  &    0     &   0   & 3&4&2&5 \\
  $345<456$: &   0   &   0   &   $-$    &   0   & 4&5&3&6 \\
  $345<347$: &   0   &   0   &    0     &  $-$  & 3&4&5&7
\end{tabular}

\medskip\noindent %
If $[257]=-$ or $[257]=0$, we can find a positive combination of the
rows of this matrix that sums to zero, regardless of the actual values
of the entries. Therefore, $[257]=+$ in any realization of~$N\!R_1^4$.
By Lemma~\ref{lem:sign}(\ref{lem:signs}), we deduce that therefore
$[157]=-$.

Now consider the rows

\medskip\qquad
\begin{tabular}[b]{c|cccc|*{4}{c@{\;}}}
             & $h_1$ & $h_2$ & $h_6$ & $h_7$
             & $i$ & $j$ & $k$ & $\ell$ \\\hline
  $127<125$: & $-[257]=-$ & $[157]=-$ & 0 & $+$  & 1&2&7&5 \\
  $145<345$: &  $+$  &    0   &    0     &   0   & 4&5&1&3 \\
  $234<345$: &   0   &   $+$  &    0     &   0   & 3&4&2&5 \\
  $345<347$: &   0   &    0   &    0     &  $-$  & 3&4&5&7
\end{tabular}\;,

\medskip\noindent
which admit a positive combination that sums to zero and therefore prove
the nonrealizability of the orientation~$N\!R_1^4$. \hfill $\Box$

\begin{remark}
  Proposition~\ref{prop:nr14} provides an example of a non-realizable
  abstract objective function that satisfies the Holt--Klee
  conditions, on a simple $4$-polytope with only $7$ facets.  The
  first examples for this were obtained on a $7$-dimensional polytope
  with $9$ facets, by G\"artner et al.\ \cite{gaertner01:_one}; Morris
  \cite{morris02:_distin} obtained examples on the $4$-cube (with $8$
  facets).  No such examples of dimension $d\le 3$ exist (Mihalisin \&
  Klee \cite{MihalisinKlee}).
\end{remark}

\begin{proposition}\label{prop:nr6}
  No Hamilton HK~AOF of $C_6(9)^\Polar$ is realizable.
\end{proposition}

\noindent\emph{Proof.} The reasoning is analogous to the proof of
Proposition~\ref{prop:nr14}; we will give the details only for
$N\!R_1^6$, and sketch the proof for the rest of the orientations.

Suppose then that we are given a realization of the polytope
$C_6(9)^\Polar$ along with a linear objective function that
induces~$N\!R_1^6$ on its graph. After an affine transformation of the
extended Gale diagram, we may suppose that $h_3=h_4=h_5=0$, where we
consider the lifting heights numbered as $h_0$,~$h_1$,\dots,$h_8$.

Now consider the rows

\medskip\qquad
\begin{tabular}[b]{c|cccccc|*{4}{c@{\;}}}
             & $h_0$ & $h_1$ & $h_2$ & $h_6$ & $h_7$ & $h_8$
             & $i$ & $j$ & $k$ & $\ell$ \\\hline
  $567<056$: & $-$ & 0 & 0 & $-[057]$ & $+$ & 0 & 5&6&7&0 \\
  $034<345$: & $+$ & 0 & 0 &    0     &  0  & 0 & 3&4&0&5 \\
  $345<456$: &  0  & 0 & 0 &   $-$    &  0  & 0 & 4&5&3&6 \\
  $345<347$: &  0  & 0 & 0 &    0     & $-$ & 0 & 3&4&5&7
\end{tabular}\;,

\medskip\noindent
from which we deduce as above that $[057]=+$, and via
Lemma~\ref{lem:sign}(\ref{lem:signs}) that $[578]=-$. But now we reach a
contradiction
via

\medskip\qquad
\begin{tabular}[b]{c|cccccc|*{4}{c@{\;}}}
             & $h_0$ & $h_1$ & $h_2$ & $h_6$ & $h_7$ & $h_8$
             & $i$ & $j$ & $k$ & $\ell$ \\\hline
  $078<058$: & $[578]=-$ & 0 & 0 & 0 & $+$ & $-[057]=-$ & 0&8&7&5 \\
  $034<345$: &    $+$    & 0 & 0 & 0 &  0  & 0          & 3&4&0&5 \\
  $345<347$: &     0     & 0 & 0 & 0 & $-$ & 0          & 3&4&5&7 \\
  $458<345$: &     0     & 0 & 0 & 0 &  0  &    $+$     & 4&5&8&3
\end{tabular}\;,

\medskip\noindent
which proves the claim. Some ``good'' sets of vanishing heights for the
remaining orientations are as follows:

\begin{center}
  \begin{tabular}{l|ccccc}
    Orientation: \rule[-1.2ex]{0pt}{0pt} & 
    $N\!R_2^6$ & $N\!R_3^6$ & $N\!R_4^6$ & $N\!R_5^6$ &
    $N\!R_6^6$ \\\hline
    Height indices: \rule{0pt}{2.2ex} & 0,1,2 & 0,5,6 & 0,5,6 & 0,1,6 & 0,1,2
  \end{tabular}
\end{center}

\medskip\noindent
This concludes the proof.\hfill $\Box$

\begin{proposition}\label{prop:nr24}
  The Hamilton HK~AOFs $N\!R_2^4$ and $N\!R_3^4$ are not
  realizable.
\end{proposition}

\noindent \emph{Proof.} The method used in the proof of
Propositions~\ref{prop:nr14} and~\ref{prop:nr6} does not directly
apply here, as no choice of vanishing heights immediately yields a
Farkas contradiction for these orientations. Therefore we prove the
nonrealizability of $N\!R_2^4$ in the following way:

Suppose that in a realization of $N\!R_2^4$, we have $[136]=+$,
and therefore $[137]=-$ by Lemma~\ref{lem:sign}(\ref{lem:signs}). This
leads to a
contradiction by the following table  for $h_1=h_4=h_5=0$:

\medskip\qquad
\begin{tabular}[b]{c|cccc}
             & $h_2$ & $h_3$ & $h_6$ & $h_7$ \\\hline
  $367<167$: & 0 & $+$ & $-[137]=+$ & $[136]=+$ \\
  $145<345$: & 0 & $-$ &     0      &    0      \\
  $145<456$: & 0 &  0  &    $-$     &    0      \\
  $145<147$: & 0 &  0  &     0      &   $-$
\end{tabular}

\medskip\noindent %
We deduce that $[136]=-$ or $[136]=0$ must hold in any realization
of~$N\!R_2^4$. But setting $h_2=h_3=h_6=0$ then yields the following
table,

\medskip\qquad
\begin{tabular}[b]{c|cccc}
             & $h_1$ & $h_4$ & $h_5$ & $h_7$ \\\hline
  $367<167$: & $-$ & 0 & 0 & $[136]$ \\
  $123<236$: & $+$ & 0 & 0 &    0      \\
  $367<236$: &  0  & 0 & 0 &   $+$
\end{tabular}\;,

\medskip\noindent and a global contradiction.

The same argument proves that $N\!R_3^4$ is nonrealizable. The only
difference between this orientation and $N\!R_2^4$ is that $345<347$
in $N\!R_3^4$, whereas $347<345$ in $N\!R_2^4$, but the proof of the
nonrealizability of~$N\!R_2^4$ did not use this inequality. \hfill
$\Box$

\begin{remark}
  The short non-realizability proofs above were found by computer,
  though they can be checked by hand.
  Propositions~\ref{prop:nr14}~and~\ref{prop:nr6} were found by trying
  to eliminate signs from all minors of~$A$ obtained by successively
  deleting triples of columns, while the proof of
  Proposition~\ref{prop:nr24} was obtained by moreover assuming
  various signs to be positive resp.\ negative. We presented instances
  of the shortest proofs found.
\end{remark}

\section{Some problems}

Our methods were successful for small dimensions and coranks, but they
do not yield (non-)\break existence statements or asymptotics for
large $d$ and $n-d$. Thus we leave the following problems open
for~now:
\begin{compactitem}[~$\bullet$~]
\item
Does $C_d(d+3)^\Polar$ have \emph{any} Hamilton~HK~AOFs for 
  even $d>6$? If not, this would give a purely combinatorial
proof that for some parameters $M(d,n)<M_{\rm ubt}(d,n)$.

What happens for odd $d\ge7$? 
%Is it true that $C_4(n)$ has no Hamilton~HK~AOFs for $n\ge8$?
%(Compare Pfeifle~\cite{Pfeifle-paths03}.)
\item
Is it true that $\ M(d,n)\ \ll\ M_{\rm ubt}(d,n)\ $
for large $n\ge d+3$, $d\ge6$? 
%One might try to start with the
%cases of small corank $n-d=3$ or of small dimension $d=6$.
\end{compactitem}
To demonstrate that the gaps in our asymptotic knowledge are substantial,
we note that in the ``diagonal'' case of $n=2d$ 
all we know is 
\[ 
2^d \ \le\ M(d,2d)\ \le\ M_{\rm ubt}(d,2d)\ \approx\ {2.6}^d.
\]
%In the case of constant dimension, we have
%\[
% ...
%\]
%where the lower bound is provided by ``deformed products''

\subsection*{Acknowledgements}
The first author is grateful for an initiation to the power of the
extended Gale diagram by Emo Welzl himself, during a stay at
ETH~Z\"urich.  Volker Kaibel has contributed to this paper through
numerous discussions and useful observations, and Christoph Eyrich
through uncommon typographical expertise.

\bibliographystyle{siam}
\bibliography{julian-bib,../bib/POLYref,../bib/POLYref2}

\end{document}